\documentclass[10pt]{article}
\usepackage{amsmath}
\usepackage{graphicx}
\usepackage{amsfonts}
\usepackage{amssymb}
\newtheorem{theorem}{Theorem}[section]

\newtheorem{corollary}[theorem]{Corollary}

\newtheorem{lemma}[theorem]{Lemma}

\newtheorem{proposition}[theorem]{Proposition}

\newenvironment{proof}[1][Proof]{\textbf{#1.} }{\ \rule{0.5em}{0.5em}}
\newdimen\dummy
\dummy=\oddsidemargin
\begin{document}

\title{Kontsevich quantization and invariant distributions on Lie groups}
\author{Martin Andler\thanks{D\'{e}partement de Mathematiques, Universit\'{e} de
Versailles Saint Quentin, \ 78035 Versailles C\'{e}dex}, Alexander
Dvorsky\thanks{Mathematics Department, Rutgers University, New Brunswick, NJ
08903} and Siddhartha Sahi\footnotemark[2]}
\maketitle
\begin{abstract}
We study Kontsevich's deformation quantization for the dual of a
finite-dimensional real Lie algebra (or superalgebra) $\frak{g}$. In this case
the Kontsevich $\star$-product defines a new convolution on $S(\frak{g})$,
regarded as the space of distributions supported at $0\in\frak{g}$. For $p\in
S(\frak{g})$, we show that the convolution operator $f\longmapsto p\star f$ is
a differential operator with \emph{analytic} germ.

We use this fact to prove a conjecture of Kashiwara and Vergne on invariant
distributions on a Lie group $G$. This  yields a new proof of Duflo's result
on local solvability of bi-invariant differential operators on a Lie group
$G$. Moreover, this new proof extends to Lie supergroups. 
\end{abstract}

\section{Introduction}

\ \ In recent work \cite{kontsevich}, M. Kontsevich has established a
remarkable result on the formality of the Hochschild complex of a smooth
manifold. An important consequence of this result is an explicit construction
of an associative $\star$-product on an arbitrary smooth Poisson manifold
$(X,\gamma)$. This $\star$-product is given by a formal power series of
bidifferential operators (depending on a parameter $\hbar$).

For a general Poisson manifold this series does not converge, except for
$\hbar=0$. Thus the analytic properties of the $\star$-product are somewhat
obscure. However, as we shall show in this paper, if $X=\frak{g}^{\ast}$ is
the dual of a finite-dimensional real Lie algebra $\frak{g}$ (or, in general,
if the coefficients of the Poisson bracket $\gamma$ are linear), then the
situation is much nicer.

In this case, we can set $\hbar=1$, and Kontsevich's $\star$-product descends
to an actual product on the algebra of polynomial functions on $\frak{g}%
^{\ast}$, which can be naturally identified with the symmetric algebra
$\mathcal{S}=\mathcal{S}(\frak{g})$. In this paper, we will identify
$\mathcal{S}(\frak{g})$ with the convolution algebra of distributions on
$\frak{g}$ supported at $0$. Then the $\star$-product can be viewed as a new
associative convolution product of such distributions.

Let $\frak{D}$ be the algebra of germs at $0$ of differential operators on
$\frak{g}$ with analytic coefficients. By virtue of the pairing between
functions and distributions, $\mathcal{S}$ is naturally a \emph{right
}$\frak{D}$-module and to emphasize this fact we shall write the action on
distributions \emph{on the right}. Our first main result is the following

\begin{theorem}
Given $p\in\mathcal{S}$ of order $l$, there is a unique element $\partial
_{p}^{\star}$ of order $l$ in $\frak{D}$ such that for all $r\in\mathcal{S} $%
\[
r\star p=r\cdot\partial_{p}^{\star}\text{ .}%
\]
\label{-main-d}
\end{theorem}

Next, in \cite[8.3.3]{kontsevich}, Kontsevich introduces a certain formal
power series
\begin{equation}
S_{1}(x)=\exp\left(  \sum_{k=1}^{\infty}c_{2k}^{(1)}\operatorname{tr}\left[
\left(  \operatorname*{ad}x\right)  ^{2k}\right]  \right)  , \label{=s1}%
\end{equation}
where the constants $c_{2k}^{(1)}$ are expressed as integrals of smooth
differential forms over certain manifolds with corners (compactified
configuration spaces, introduced in \cite{fulton-macph}).

\begin{proposition}
The power series $S_{1}(x)$ converges to an analytic function $\tau(x)$ in
some neighborhood of $0$ in $\frak{g}$. \label{-tau}
\end{proposition}

For $a\in\frak{g}$, define the adjoint vector field $\operatorname*{adj}{}%
_{a}$ by the formula
\[
\operatorname*{adj}{}_{a}(f)(x)=\frac{d}{dt}f(\exp(-ta)\cdot x)|_{t=0}\,.
\]
Let
\[
\mathcal{I}=\mathcal{S}^{\frak{g}}=\left\{  p\in\mathcal{S\,}\left|
\,p\cdot\operatorname*{adj}{}_{a}=0\text{ for all }a\in\frak{g}\right.
\right\}
\]
be the subalgebra of invariant distributions, and let $\frak{R}$ be the right
ideal of $\frak{D}$ generated by the germs of the differential operators
$\operatorname*{adj}{}_{a}$, $a\in\frak{g}$.

Our next result is the following

\begin{theorem}
For $p\in\mathcal{I}$, the operator
\[
T_{p}=\partial_{p}\tau-\tau\partial_{p\tau}^{\star}%
\]
belongs to $\frak{R}$. \label{-main-r}
\end{theorem}

Here $p\tau\in\mathcal{S}$ and $\partial_{p\tau}^{\star}$ is the corresponding
differential operator defined in Theorem \ref{-main-d}, and for $r\in
\mathcal{S\;\;}r\cdot\partial_{p}=r\ast_{\frak{g}}p$ (convolution of
distributions on $\frak{g}$).

We now discuss some applications of these results to invariant distributions
on Lie groups. Let $G$ be a finite-dimensional Lie group with Lie algebra
$\frak{g}$. We write $\mathcal{U}=\mathcal{U}(\frak{g})$ for the enveloping
algebra of $\frak{g,}$ and $\mathcal{Z}=\mathcal{U}^{\frak{g}}$ for the center
of the enveloping algebra. We identify $\mathcal{U}$ with the convolution
algebra of distributions on $G$ supported at $1\in G$.

Let $\mathbf{U}$ and $\mathbf{S}$ be the spaces of germs at $1$ and $0$ of
distributions on $G$ and $\frak{g}$, and let $\mathbf{Z}$ and $\mathbf{I}$ be
the subspaces of $\frak{g}$-invariant germs. We denote by $\exp_{\ast
}:\mathbf{S\rightarrow U}$ and $\exp^{\ast}:\mathbf{U\rightarrow S}$ the
pushforward and pullback of distributions (germs) under the exponential map;
and define $\eta:\mathbf{S\rightarrow U}$ by
\begin{equation}
\eta(p)=\exp_{\ast}(pq)\text{,} \label{=duflo-eta}%
\end{equation}
where for $x\in\frak{g}$%

\[
q(x):=\det\left(  \frac{e^{\operatorname{ad}(x)/2}-e^{-\operatorname{ad}%
(x)/2}}{\operatorname{ad}(x)}\right)  ^{1/2}=\exp\left(  \sum_{k=1}^{\infty
}\frac{B_{2k}}{4k(2k)!}\operatorname{tr}\left[  \left(  \operatorname*{ad}%
x\right)  ^{2k}\right]  \right)  \text{.}%
\]
Here $B_{2k}$ are the Bernoulli numbers.

It follows from \cite[Ch. 8]{kontsevich}, that the restriction of $\eta$ to
$\mathcal{I}$ is an \emph{algebra} isomorphism from $\mathcal{I}$ to
$\mathcal{Z}$, i.e.,
\begin{equation}
\eta(p_{1}\ast_{\frak{g}}p_{2})=\eta(p_{1})\ast_{G}\eta(p_{2})
\label{-duflo-eta}%
\end{equation}
for $p_{1},\,p_{2}\in\mathcal{I}$, where $\ast_{\frak{g}}$ and $\ast_{G}$
denote convolutions in the Lie algebra and Lie group respectively. (This fact
was first established by Duflo in \cite{duflo-iso}$^{^{?}}$ using Kirillov's
orbit method).

From Theorems \ref{-main-d} and \ref{-main-r}\thinspace, we can deduce

\begin{theorem}
For $p\in\mathcal{I}$, the following differential operator $D_{p}%
:\mathbf{S\rightarrow S}$ lies in $\frak{R\,}$:
\[
P\cdot D_{p}=\exp^{\ast}(\eta(P\ast_{\frak{g}}p)-\eta(P)\ast_{G}\eta(p)).
\]
\label{-kv-conjecture}
\end{theorem}

This result was conjectured by M. Kashiwara and M. Vergne in \cite{kash-ver},
and proved by them for $\frak{g}$ solvable. It has several important
corollaries. First, since operators from $\frak{R}$ annihilate distributions
from $\mathbf{I}$, we obtain the following extension of the isomorphism
(\ref{=duflo-eta}):

\begin{theorem}
For $p\in\mathcal{I}$ \ and $P\in$ $\mathbf{I}$, we have $\eta(P\ast
_{\frak{g}}p)=\eta(P)\ast_{G}\eta(p)$. \endproof\label{-module-iso}
\end{theorem}

Next, recall that a distribution $Q$ on $G$ is called an eigendistribution if
there is a character $\chi:\mathcal{Z}\longrightarrow\mathbb{C}$ such that
\[
Q\ast_{G}z=\chi(z)Q\text{ for all }z\text{ in }\mathcal{Z}\text{.}%
\]
Eigendistributions on the Lie algebra are defined similarly, and we obtain the following

\begin{theorem}
The map $\eta$ takes germs of invariant eigendistributions on the Lie algebra
to those on the Lie group.\endproof\label{-eigendist}
\end{theorem}

As another consequence of Theorem \ref{-module-iso}, we obtain

\begin{theorem}
Every nonzero bi-invariant differential operator on $G$ admits a (local)
fundamental solution.\label{-fund-sol}
\end{theorem}

Note that the analogous statement for left-invariant differential operators on
$G$ is known to be false. In fact, the first example of a not locally solvable
partial differential operator (\cite{lewy}) is a left-invariant vector field
on a three-dimensional Heisenberg group. 

\begin{corollary}
Every nonzero bi-invariant differential operator on $G$ is locally solvable. \endproof
\end{corollary}

This last result was originally proved by Duflo in \cite{duflo-inv} (see also
\cite{helgason}, \cite{rais-nilpotent}, \cite{duflo-rais} for various special
cases). An interpretation of the statement of Theorem \ref{-fund-sol} in terms
of Kirillov's orbit method can be found in \cite[6.1 and 8.2]{kirillov}.

Most of our results extend to the case of finite-dimensional Lie superalgebras
without much additional effort. For simplicity of exposition, we work with
ordinary Lie algebras and at the end explain the modifications necessary for
the super-setting.\smallskip

\textbf{Acknowledgments. --} We would like to thank I. Gelfand, F. Knop and
the other participants of the Gelfand seminar at Rutgers in which we first
learned the main ideas behind the Kontsevich formality theorem. The authors
would also like to thank M. Kontsevich for sharing his insights with us and P.
Deligne for many valuable suggestions.\label{-sec-intro}

\section{Relevant facts from \cite{kontsevich}}

\subsection{Configuration spaces}

Following Kontsevich, for $2n+m\geq2$, we define
\[
\operatorname{Conf}_{n,m}=\{(p_{1},\ldots,p_{n};q_{1},\ldots,q_{m}%
)|\,\,p_{i}\in\mathcal{H},p_{i_{1}}\neq p_{i_{2}},q_{j}\in\mathbb{R},q_{j_{1}%
}\neq q_{j_{2}}\}.
\]
Here $\mathcal{H}=\{z\in\mathbb{C}:\operatorname{Im}z>0\}$.

$\operatorname{Conf}_{n,m}$ is a smooth manifold of dimension $2n+m$, with a
free action of the two-dimensional group $G^{1}=\{z\rightarrow
az+b\,|\,a>0,b\in\mathbb{R}\}.$ We consider the quotient space
\[
C_{n,m}=\operatorname{Conf}_{n,m}/G^{1}.
\]
Then $\dim C_{n,m}=2n+m-2$.

The spaces $\operatorname{Conf}_{n,m}$ and $C_{n,m}$ have $m!$ connected
components, and we denote by $\operatorname{Conf}_{n,m}^{+}$ and $C_{n,m}^{+}$
the component where
\[
q_{1}<q_{2}<\ldots<q_{m}\text{.}%
\]

Similarly, for $n\geq2$, we define
\begin{align*}
\operatorname{Conf}_{n}  &  =\{(p_{1},\ldots,p_{n})\text{ }|\,\,p_{i}%
\in\mathbb{C}\text{, }p_{i_{1}}\neq p_{i_{2}}\}\text{ and}\\
C_{n}  &  =\operatorname{Conf}_{n}/G^{2}\text{, where }G^{2}=\{z\rightarrow
az+b\,|a>0,b\in\mathbb{C}\}.
\end{align*}

Clearly, $\dim C_{n}=2n-3.$

Consider, for example, the space
\[
C_{2,0}=\{(p_{1},p_{2})\in\mathcal{H}^{2}\,|p_{1}\neq p_{2}\}/G^{1}\text{.}%
\]
For each point $c\in C_{2,0}$ we \ can choose a unique representative of the
form $(\sqrt{-1},z)\in\operatorname{Conf}_{2,0}$. Thus $C_{2,0}$ is
homeomorphic to $\mathcal{H}\setminus\{\sqrt{-1}\}$.

Similarly, it is easy to see that
\begin{align*}
C_{2}  &  \approx S^{1},\\
C_{1,1}  &  \approx(0,1),\\
C_{0,2}  &  \approx\{0,1\}\text{.}%
\end{align*}

The spaces $C_{n,m}$ have natural compactifications $\overline{C}_{n,m}$, such
that the boundary components in $\overline{C}_{n,m}\setminus C_{n,m}$
correspond to various degenerations of the configurations of points.

For example,
\begin{align*}
\overline{C}_{2}  &  =C_{2},\\
\overline{C}_{1,1}  &  =C_{1,1}\sqcup C_{0,2}=[0,1],
\end{align*}

Of particular interest is the space
\[
\overline{C}_{2,0}=C_{2,0}\sqcup(C_{0,2}\sqcup C_{1,1}\sqcup C_{1,1}\sqcup
C_{2}),
\]
which can be drawn as ``\textsc{the Eye}''.%

\begin{figure}
[ptbh]
\begin{center}
\includegraphics[
trim=0.000000in 0.330945in 0.000000in 0.000000in,
height=0.8216in,
width=1.7928in
]%
{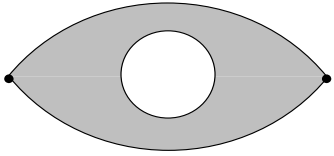}%
\caption{$\overline{C}_{2,0}$ (\textsc{the Eye})}%
\end{center}
\end{figure}

The circle $C_{2}$ represents two points coming close together in the interior
of $\mathcal{H}$, and two arcs $C_{1,1}\sqcup C_{1,1}$ represent the first
point (or the second point) coming close to the real line. Finally, the two
corners $C_{0,2}$ correspond to both points approaching the real line.

The compactification $\overline{C}_{n,m}$ is described in general in \cite[Ch.
5]{kontsevich}, and the boundary components are products of various $C_{k,l}%
$'s and $C_{p}$'s$.$

In particular, the boundary components of codimension 1 in $\overline{C}%
_{n,m}$ are of the form

\begin{enumerate}
\item $C_{k}\times C_{n-k+1,m}$ ($k\geq2$) or

\item $C_{k,l}\times C_{n-k,m-l+1}$ ($k\geq1$ or $l\geq2$).
\end{enumerate}

The first case corresponds to a cluster of $k$ points in $\mathcal{H}$ coming
infinitely close to each other in $\mathcal{H}$ , and the second case
corresponds to a cluster of $k$ points in $\mathcal{H}$ and $l$ points in
$\mathbb{R}$ coming infinitely close to each other and $\mathbb{R}$.

\subsection{Graphs and weights}

We consider the class of ``admissible'' oriented graphs $\left\{
\Gamma=(V_{\Gamma},E_{\Gamma})\right\}  $ such that

\begin{itemize}
\item $V_{\Gamma}$ consists of $n$ vertices of the ``first'' type labeled
$\{1,\ldots,n\}\ $and $m$ vertices of the ``second'' type labeled
\{$\overline{1},\overline{2},\ldots,\overline{m}\}$.

\item  For each vertex $i$ of the first type, there are $k_{i}$ edges in
$E_{\Gamma}$ starting from $i$, labeled
\[
\{e_{i}^{1},e_{i}^{2},...,e_{i}^{k_{i}}\}.
\]

\item $E_{\Gamma}$ contains no edges starting from the vertices of second
type, and no loops (edges of type $(v,v)$) or multiple edges.
\end{itemize}

For each point in $(p_{1},\ldots,p_{n};q_{1},\ldots,q_{m})\in
\operatorname{Conf}_{n,m}$, we can draw the graph $\Gamma$ in the closed upper
half-plane $\overline{\mathcal{H}}=\mathcal{H}\sqcup\mathbb{R}$ as follows:

\begin{quote}
vertices $i$ are placed at $p_{i}\in\mathcal{H}$, vertices $\overline{j}$ are
placed at $q_{j}\in$ $\mathbb{R}$, and edges are drawn as Lobachevsky
geodesics (semicircles centered on $\mathbb{R}$) connecting the vertices.
\end{quote}

For each edge $e\in E_{\Gamma}$ connecting the vertex at $z_{1}\in\mathcal{H}$
to $z_{2}\in\overline{\mathcal{H}}$, we define $\phi_{e}$ to be the angle
between $e$ and the vertical line through $z_{1}$. Then
\begin{equation}
\phi_{e}=\frac{1}{2\sqrt{-1}}\operatorname{Log}\frac{(z_{2}-z_{1}%
)(\overline{z}_{2}-z_{1})}{(z_{2}-\overline{z}_{1})(\overline{z}_{2}%
-\overline{z}_{1})}\text{.} \label{=angle}%
\end{equation}

Define a map $\Phi_{\Gamma}$ from $\operatorname{Conf}_{n,m}\ $to the
$k$-dimensional torus $\mathbb{T}^{k}$, where $k=\sum_{i=1}^{n}k_{i}$ as
follows
\begin{equation}
\Phi_{\Gamma}(p_{1},\ldots,p_{n};q_{1},\ldots,q_{m})=(\phi_{e_{1}^{1}}%
,\phi_{e_{1}^{2}},\ldots,\phi_{e_{n}^{k_{n}}})\text{.} \label{=map-to-t}%
\end{equation}
Observe that $\Phi_{\Gamma}$ descends to a map from $C_{n,m}$ to
$\mathbb{T}^{k}$. Define $\omega_{\Gamma}$ to be the pullback to $C_{n,m}$ of
the volume form on $\mathbb{T}^{k}$.

The form $\omega_{\Gamma}$ is a smooth $k$-form on $C_{n,m}$, which extends
continuously to $\overline{C}_{n,m}$. We define the \emph{weight}%
\footnote{This weight differs from $w_{\Gamma}$ on p.6 of \cite{kontsevich} by
the factor of $\frac{1}{n!}$ and from the weight $W_{\Gamma}$ on p. 23 by a
factor of $\frac{1}{k_{1}!\ldots k_{n}!}$.}\emph{ }$w_{\Gamma}$ by
\[
w_{\Gamma}=\int\limits_{\overline{C}_{n,m}^{+}}\omega_{\Gamma}\text{.}%
\]
Evidently, $w_{\Gamma}=0$, unless $k=\dim C_{n,m}=2n+m-2$.\label{ss-weights}

\subsection{Star-product}

Let
\[
\gamma(f_{1},f_{2})=\sum_{i,j=1}^{d}\gamma^{ij}\partial_{i}(f_{1})\partial
_{j}(f_{2})
\]
be a smooth Poisson bracket (bivector field) on $C^{\infty}(\mathbb{R}^{d})$.
Kontsevich's main result gives the following construction of an
\emph{associative star-product} on $C^{\infty}(\mathbb{R}^{d})\left[  \left[
h\right]  \right]  $:
\begin{equation}
f_{1}\star f_{2}=\sum_{n=0}^{\infty}\frac{\hbar^{n}}{n!}B_{n}(f_{1}%
,f_{2})\text{,} \label{=star-pr}%
\end{equation}
where $B_{0}(f_{1},f_{2})=f_{1}f_{2}$, $B_{1}(f_{1},f_{2})=\gamma(f_{1}%
,f_{2})$ and in general
\begin{equation}
B_{n}(f_{1},f_{2})=\sum_{\Gamma\in G_{n}}w_{\Gamma}B_{\Gamma}(f_{1}%
,f_{2})\text{.} \label{=star-summation}%
\end{equation}
Here $G_{n}$ is defined as a class of admissible graphs with $n$ vertices of
the first type, two vertices of the second type and $k_{i}=2,$ $1\leq i\leq
n$. The weight $w_{\Gamma}$ is as in the previous subsection, and $B_{\Gamma}$
is a bidifferential operator defined as follows:
\begin{align}
B_{\Gamma}(f_{1},f_{2})  &  =\sum_{I:E_{\Gamma}\rightarrow\{1,2,\ldots
,d\}}\left[  \prod_{k=1}^{n}\left(  \prod_{e\in E_{\Gamma},e=(\ast,k)}%
\partial_{I(e)}\right)  \gamma^{I(e_{k}^{1})I(e_{k}^{2})}\right]
\times\label{=ouch}\\
&  \times\left(  \prod_{e\in E_{\Gamma},e=(\ast,\overline{1})}\partial
_{I(e)}\right)  f_{1}\times\left(  \prod_{e\in E_{\Gamma},e=(\ast,\overline
{2})}\partial_{I(e)}\right)  f_{2}.\nonumber
\end{align}
To obtain this operator, we relabel each vertex $k,$ $1\leq k\leq n$ by a
component $\gamma^{i_{k}j_{k}}$ of $\gamma$, relabel $\overline{1}$ as $f_{1}$
and $\overline{2}$ as $f_{2}$. We also label the edge $e_{k}^{1}$ with $i_{k}$
and the edge $e_{k}^{2}$ with $j_{k}$. Then an edge leading to the vertex
indicates differentiation of the function labeling the vertex with respect to
the index labeling the edge.

Below is an example of a graph from $\Gamma_{0}\in G_{2}$ with
\[
e_{1}^{1}=(1,\overline{1}),\,e_{1}^{2}=(1,2),\,e_{2}^{1}=(2,\overline
{1}),e_{2}^{2}=(2,\overline{2}).
\]%
\[
\text{\ }%
\begin{array}
[c]{ccc}%
\overset{1}{\bullet} & \rightarrow & \overset{2}{\bullet}\\
\downarrow & \swarrow & \downarrow\\
\underset{\overline{1}}{\bullet} &  & \underset{\overline{2}}{\bullet}%
\end{array}
\text{\ or (relabeling)\ \ \ }%
\begin{array}
[c]{ccc}%
\overset{\gamma^{i_{1}j_{1}}}{\bullet} & \overset{_{j_{1}}}{\rightarrow} &
\overset{\gamma^{i_{2}j_{2}}}{\bullet}\\
\,\downarrow^{_{i_{1}}} & \overset{i_{2}}{\swarrow} & \downarrow^{j_{2}}\\
\underset{f_{1}}{\bullet} &  & \underset{f_{2}}{\bullet}%
\end{array}
.
\]

The corresponding bidifferential operator $B_{\Gamma_{0}}$ is given by
\begin{equation}
B_{\Gamma_{0}}(f_{1},f_{2})=\sum_{1\leq i_{1},j_{1},i_{2},j_{2}\leq d}%
\gamma^{i_{1}j_{1}}\partial_{j_{1}}(\gamma^{i_{2}j_{2}})\,\,\partial_{i_{1}%
}\partial_{i_{2}}(f_{1})\,\,\partial_{j_{2}}(f_{2})\text{.}
\label{=ouch-example}%
\end{equation}

\subsection{Tangent map}

Let $\gamma=\sum_{i,j=1}^{d}\gamma^{ij}\partial_{i}\wedge\partial_{j}$ be a
Poisson bivector field on $\mathbb{R}^{d}$ as before; following \cite[Ch.
8]{kontsevich}, we define the map $I_{T}$ from $C^{\infty}(\mathbb{R}^{d})$ to
$C^{\infty}(\mathbb{R}^{d})\left[  \left[  h\right]  \right]  $ by
\[
I_{T}(f)=\sum_{n=0}^{\infty}\frac{\hbar^{n}}{n!}\sum_{\Gamma\in H_{n}%
}w_{\Gamma}D_{\Gamma}(f)\text{.}%
\]
Here $H_{n\text{ }}$ is the set of all admissible graphs with $n+1$ vertices
of first type and no vertices of the second type such that $k_{i}=2,$ $1\leq
i\leq n$ and $k_{n+1}=0$. The differential operator $D_{\Gamma}$ is defined by
labeling the vertex $n+1$ with $f$, the remaining vertices by components of
$\gamma$, and then proceeding as in formula (\ref{=ouch}).

Let $\mathcal{I}$ be the center of the Poisson algebra $(C^{\infty}%
(\mathbb{R}^{d}),\gamma)$, i.e.
\[
\mathcal{I}=\left\{  f\in C^{\infty}(\mathbb{R}^{d}):\gamma(f,g)=0\text{ for
all }g\in C^{\infty}(\mathbb{R}^{d})\right\}  \text{.}%
\]
Then, as established in \cite[Ch. 8]{kontsevich},
\begin{equation}
I_{T}\left(  f_{1})\star I_{T}(f_{2}\right)  =I_{T}(f_{1}f_{2})
\label{=cup-on-functions}%
\end{equation}
for $f_{1},f_{2}\in\mathcal{I}$.

Let us briefly explain the ideas underlying the above equality.

Kontsevich's main result is a construction of an $L^{\infty}$-quasiisomorphism
$\mathcal{U}$ between the differential graded Lie algebras $\mathcal{T}$ and
$\mathcal{D}$ consisting of the skew polyvector fields and polydifferential
operators on $\mathbb{R}^{d}$, respectively. The quasiisomorphism
$\mathcal{U}$ can be regarded as a (geometric) map between formal manifolds
($Q$-manifolds) associated to $\mathcal{T}$ \-\thinspace and $\mathcal{D}$, respectively.

The skew bivector field $\hbar\gamma$ can be regarded as a point distribution
$\delta_{\hbar\gamma}$ in the first $Q$-manifold. Then the tangent spaces at
$\hbar\gamma$ and $\mathcal{U}(\hbar\gamma)$ are differential graded Lie
algebras $\mathcal{T}^{\prime}\,$and $\mathcal{D}^{\prime}$, with the same
underlying Lie algebras as $\mathcal{T}$ and $\mathcal{D}$, but with modified
differentials. The derivative $d\mathcal{U}$ of $\mathcal{U}$ at $\hbar\gamma$
induces a graded linear map between $\mathcal{T}^{\prime}$ and $\mathcal{D}%
^{\prime}$, and $I_{T}$ is precisely the restriction of $d\mathcal{U}$ to the
$0$-vector fields (i.e., functions).

The equality (\ref{=cup-on-functions}) follows from compatibility of
$d\mathcal{U}$ with the natural cup product structures on the (tangent)
cohomology of $\mathcal{T}^{\prime}$ and $\mathcal{D}^{\prime}$, as
established in \cite[Theorem 8.2]{kontsevich}.

The argument proving \cite[Theorem 8.2]{kontsevich}, which was merely outlined
in \cite{kontsevich}, was communicated to us by M. Kontsevich
(\cite{maxims-note}), and actually proves a slightly stronger result. Since we
need this stronger version, we shall reproduce Kontsevich's argument for the
case of Lie algebras in Section \ref{-tang-for-lie}.

\section{Star-product for Lie algebras}

In the rest of the paper we restrict ourselves to the situation when $X=$
$\frak{g}^{\ast}$ is the dual of a finite dimensional Lie algebra $\frak{g}$,
and $\gamma$ is the standard Poisson bracket \ on $\frak{g}^{\ast}$. In this
case the coefficients of a bivector field $\gamma\,$are linear functions on
$\frak{g}^{\ast}$. More precisely, let $x_{1},x_{2},\ldots,x_{d}$ be some
basis of $\frak{g}$ viewed as coordinate functions on $\frak{g}^{\ast}$, and
let $\{c_{ij}^{k}\}$ be the corresponding set of structure constants.

We normalize $x_{1},x_{2},\ldots,x_{d}$ such that
\begin{equation}
\left|  c_{ij}^{k}\right|  \leq2\text{ for all }i,j\text{ and }k.
\label{=normalize}%
\end{equation}
Then
\[
\left[  x_{i},x_{j}\right]  =\sum_{k=1}^{d}c_{ij}^{k}x_{k}\text{,}%
\]
and
\[
\gamma(f_{1},f_{2})=\frac{1}{2}\sum_{i,j,k}c_{ij}^{k}x_{k}\frac{\partial
f_{1}}{\partial x_{i}}\frac{\partial f_{2}}{\partial x_{j}}\text{,}%
\]
i.e.,
\[
\gamma^{ij}=\frac{1}{2}\sum_{k=1}^{d}c_{ij}^{k}x_{k}\text{.}%
\]

The linearity of the coefficients of $\gamma$ implies that the formula
(\ref{=star-summation}) for the bidifferential operator $B_{n}$ can be
rewritten as
\[
B_{n}=\sum_{\Gamma\in A_{n}}w_{\Gamma}B_{\Gamma}\text{,}%
\]
where $A_{n}$ consists of those graphs in $G_{n}$ for which there is at most
one incoming edge at every vertex of the first type. Indeed, if $\Gamma\in
G_{n}\setminus A_{n}$, the corresponding bidifferential operator
$B_{\Gamma,\gamma}$ is automatically $0$.

Next, let $f_{1}$ and $f_{2}$ be two polynomials on $\frak{g}^{\ast}$ with
$\deg(f_{1})=l_{1}$, $\deg(f_{2})=l_{2}$. Then we remark that the graphs
contributing to the star-product formula for $f_{1}\star f_{2}$ can have no
more than $l_{1}+l_{2}$ vertices of the first type. Indeed, for any $\Gamma\in
A_{n}$ the corresponding bidifferential operator $B_{\Gamma,\gamma}$ contains
exactly $2n$ differentiations. When $2n>n+l_{1}+l_{2}$, $B_{\Gamma}$ is
obviously $0$ (because $f_{1}$ can be differentiated at most $l_{1}$ times,
$f_{2}$ at most $l_{2}$ times and each of the coefficients of $\gamma$
corresponding to the remaining vertices at most once). Hence
\[
f_{1}\star f_{2}=\sum_{n=0}^{l_{1}+l_{2}}\frac{\hbar^{n}}{n!}\sum_{\Gamma\in
A_{n}}w_{\Gamma}B_{\Gamma}(f_{1},f_{2})\text{.}%
\]

This sum is finite, and if we set $\hbar=1$, we obtain a polynomial on
$\frak{g}^{\ast}$ of degree $l_{1}+l_{2}$. Therefore, Kontsevich's $\star
$-product descends to an \emph{actual product} on the algebra of polynomial
functions on $\frak{g}^{\ast}$, which can be naturally identified with the
symmetric algebra $\mathcal{S}=\mathcal{S}(\mathcal{\frak{g}})$.

Let $\frak{W}\left(  \frak{g}\right)  $ and $\frak{W}\left(  \frak{g}^{\ast
}\right)  $ be the Weyl algebras of polynomial coefficient differential
operators on $\frak{g}$ and $\frak{g}^{\ast}$ respectively. Both of these
algebras are generated by $\frak{g}\oplus\frak{g}^{\ast}$, and there is a
unique anti-isomorphism $\iota$ between them, which is the identity on the
generators. Observe that $\iota$ interchanges multiplication and
differentiation operators.

We will now change our point of view and regard $\mathcal{S}$ as the algebra
of distributions with point support (at $0\in\mathcal{\frak{g}}$). Then
$\mathcal{S}$ is naturally a \emph{right module} for the algebra
$\frak{W}\left(  \frak{g}\right)  $ via $\iota$. As discussed in the
introduction, we shall write differential operators acting on distributions,
on the right.

This approach leads to the following

\begin{lemma}
Let $p\in\mathcal{S}$ $\ $be a homogeneous distribution of order $l$, and let
$\Gamma\in A_{n}$. Then there is a differential operator $\partial_{\Gamma
}^{p}\in$ $\frak{W}\left(  \frak{g}\right)  $ of order at most $l,$ with
polynomial coefficients of degree at most $n$, such that for any
$r\in\mathcal{S}$
\[
B_{\Gamma}(r,p)=r\cdot\partial_{\Gamma}^{p}\text{.}%
\]
\label{-left-to-right}
\end{lemma}

\begin{proof}
From the formula (\ref{=ouch}) for $B_{\Gamma}(r,p)$, it is clear that
$B_{\Gamma}(\cdot,p)$ is a differential operator from $\frak{W}\left(
\frak{g}^{\ast}\right)  $. Since at most $n$ edges of $\Gamma$ can terminate
at the vertex $\overline{1}$ (corresponding to $r$), the order of $B_{\Gamma
}(\cdot,p)$ is at most $n.$ Also, the degree of the coefficients of this
differential operator is
\[
n+l-\text{\#edges terminating at vertices other than }\overline{1}\text{,}%
\]
which is at most $l$.

Now let $\partial_{\Gamma}^{p}=\iota\left(  B_{\Gamma}(\cdot,p)\right)  $.
Then the statement follows from the discussion above.
\end{proof}

\textbf{Remark. }Let $y_{1},\ldots,y_{d}$ be the basis of $\frak{g}^{\ast}$
dual to $x_{1},\ldots,x_{d},$ viewed as coordinate functions on $\frak{g}$.
Then, by the lemma above, for $p$ homogeneous of order $l$, we can write%
\begin{equation}
r\cdot\partial_{\Gamma}^{p}=r\cdot\left(  \sum_{\left|  \beta\right|  \leq
l}\left(  \sum_{\alpha}c_{\alpha\beta}^{\Gamma}y^{\alpha}\right)  \partial
_{y}^{\beta}\right)  \text{,} \label{=c-alpha-beta-Gamma}%
\end{equation}
where $\alpha$ and $\beta$ are multi-indices, and, as usual
\begin{align*}
y^{\alpha}  &  =y_{1}^{\alpha_{1}}y_{2}^{\alpha_{2}}\ldots y_{d}^{\alpha_{d}%
}\text{,}\\
\partial_{y}^{\beta}  &  =\left(  \frac{\partial}{\partial y_{1}}\right)
^{\beta_{1}}\left(  \frac{\partial}{\partial y_{2}}\right)  ^{\beta_{2}}%
\cdots\left(  \frac{\partial}{\partial y_{d}}\right)  ^{\beta_{d}}\text{.}%
\end{align*}
Moreover,
\begin{align*}
\left|  \alpha\right|   &  =\text{\#edges terminating at }\overline{1}%
\text{,}\\
\left|  \beta\right|   &  =n+l-\text{\#edges teminating at other vertices.}%
\end{align*}
Subtracting second equation from the first, obtain $\left|  \alpha\right|
-\left|  \beta\right|  =n-l$.

We now proceed with the proof of Theorem \ref{-main-d}.

\begin{lemma}
$\left|  A_{n}\right|  <(8e)^{n}n!$.\label{-count-graph}
\end{lemma}

\begin{proof}
To describe any graph from $A_{n}$, it suffices to provide the following data
for each vertex $j,$ $1\leq j\leq n$:

\begin{itemize}
\item  the source of an incoming edge, if any ($n$ choices, counting the case
when there is no incoming edge);

\item  whether the vertex is connected to the vertex $\overline{1}$, vertex
$\overline{2}$, both or neither (4 choices).
\end{itemize}

Hence we have at most $4n$ choices for each of $n$ vertices of the first type.
Taking into account the possible labelings of the outgoing edges at each
vertex, we get
\[
\left|  A_{n}\right|  \leq2^{n}(4n)^{n}\text{.}%
\]
Now the statement follows from the inequality $n^{n}<e^{n}n!$ .
\end{proof}

\begin{lemma}
Let $\Gamma\in A_{n}$. Then
\[
\left|  w_{\Gamma}\right|  \leq4^{n}\text{.}%
\]
\label{-estimate-weight}
\end{lemma}

\begin{proof}
Recall that the weight $w_{\Gamma}$ associated to the graph $\Gamma\in A_{n}$
is given by the formula
\[
w_{\Gamma}=\int_{C_{n,2}^{+}}\omega_{\Gamma}\text{,}%
\]
where $\omega_{\Gamma}$ is the pullback of the volume form from the
$2n$-dimensional torus $\mathbb{T}^{2n}$ under $\Phi_{\Gamma}$.

We consider the preimage of a generic point $\mathbf{\varphi\in}%
\mathbb{T}^{2n}$ under $\Phi_{\Gamma}$. Using the action of $G^{1}$, we can
fix $q_{1}=0$ and $q_{2}=1$, and thus identify $C_{n,2}^{+}$ with
$\mathcal{H}^{n}$.

Rewriting the formula (\ref{=angle}) as
\[
(z_{2}-z_{1})(\overline{z}_{2}-z_{1})=e^{2\sqrt{-1}\phi_{e}}(z_{2}%
-\overline{z}_{1})(\overline{z}_{2}-\overline{z}_{1})\text{,}%
\]
we see that each edge gives rise to a pair of quadratic equations involving
the coordinates of its endpoints. Thus we obtain a system of $2n$ quadratic
equations in $2n$ real variables (the real and imaginary parts of $p_{1}%
,p_{2},\ldots,p_{n}$). Generically, this system has at most $2^{2n}=4^{n}$
(complex) solutions.

Therefore, we conclude that
\[
\left|  w_{\Gamma}\right|  \leq4^{n}V\text{ ,}%
\]
where $V$ is the volume of $\Phi_{\Gamma}(C_{n,2}^{+})$ $\subset$
$\mathbb{T}^{2n}$. Obviously, $V\leq1$ and the statement follows.\smallskip
\end{proof}

\begin{proof}
[Proof of Theorem \ref{-main-d}]Without loss of generality, we may take
\[
p=x_{1}^{a_{1}}x_{2}^{a_{2}}\ldots x_{d}^{a_{d}}\text{, with }a_{1}%
+\ldots+a_{d}=l\text{.}%
\]

According to the formulas (\ref{=star-pr}), (\ref{=star-summation}) and Lemma
\ref{-left-to-right},
\begin{equation}
r\star p=r\cdot\left(
{\displaystyle\sum\limits_{n=0}^{\infty}}
\frac{1}{n!}%
{\displaystyle\sum\limits_{\Gamma\in A_{n}}^{\infty}}
w_{\Gamma}\partial_{\Gamma}^{p}\right)  \text{.} \label{=bigsum}%
\end{equation}
The expression above can be rewritten as
\[
r\star p=r\cdot\left(  \sum_{\left|  \beta\right|  \leq l}\left(  \sum
_{\alpha}c_{\alpha\beta}y^{\alpha}\right)  \partial_{y}^{\beta}\right)
\text{.}%
\]

Set $m=l+\left|  \alpha\right|  -\left|  \beta\right|  $. It follows from the
discussion after Lemma \ref{-left-to-right}, that only the graphs $\Gamma\in
A_{m}$ can contribute to $c_{\alpha\beta}$. Observe that all partial
derivatives of $\gamma^{ij}$ are at most 1 in absolute value (this follows
form the inequality (\ref{=normalize})), and that any derivative of $p$ is a
monomial with the coefficient not exceeding $C_{p}=a_{1}!a_{2}!\ldots a_{d}!$.
Therefore the absolute value of the coefficient $c_{\alpha\beta}^{\Gamma}$ in
(\ref{=c-alpha-beta-Gamma}) does not exceed $C_{p}$. 

Using Lemma \ref{-estimate-weight}, we obtain
\[
\left|  c_{\alpha\beta}\right|  =\frac{1}{m!}\left|  \sum_{\Gamma\in A_{m}%
}w_{\Gamma}c_{\alpha\beta}^{\Gamma}\right|  \leq\frac{\left|  A_{m}\right|
}{m!}4^{m}C_{p}\text{.}%
\]
Now by Lemma \ref{-count-graph}, $\left|  A_{m}\right|  <(8e)^{m}m!$, hence
\[
\left|  c_{\alpha\beta}\right|  <C_{p}(32e)^{m}\leq C_{p}^{\prime
}(32e)^{\left|  \alpha\right|  }\text{,}%
\]
where $C_{p}^{\prime}=(32e)^{l}C_{p}$.

Therefore, the series $\sum_{\alpha}c_{\alpha\beta}y^{\alpha}$ converges
absolutely in the polydisk of radius $\dfrac{1}{32e}$. Hence in this
neighborhood the formal series (\ref{=bigsum}) defines a differential operator
of the form $\sum_{\left|  \beta\right|  \leq l}M_{\beta}\partial^{\beta}$,
where $M_{\beta}$ is the operator of multiplication by some \emph{analytic
}function. The statement follows.\smallskip\label{-star-for-lie}
\end{proof}

\section{Tangent map for Lie algebras}

We start by establishing Proposition \ref{-tau}.

\begin{proof}
[Proof of Proposition \ref{-tau}]It suffices to show that the constants
$c_{2k}^{(1)\text{ }}$ from \cite[8.3.3]{kontsevich} satisfy estimates of the
form
\[
c_{2k}^{(1)}\leq c^{k}%
\]
for some constant $c$.

According to \cite[8.3.3.1]{kontsevich}, we have
\[
c_{n}^{(1)}=\frac{w_{\operatorname{Wh}_{n}}}{2^{n}},
\]
where ``the wheel'' $\operatorname{Wh}_{n}$ is the graph from $H_{n}$ shown in
Figure \ref{-fig-wheel}.%

\begin{figure}
[ptbh]
\begin{center}
\includegraphics[
trim=-0.011608in 0.011512in -0.013106in -0.011512in,
height=1.4451in,
width=1.3361in
]%
{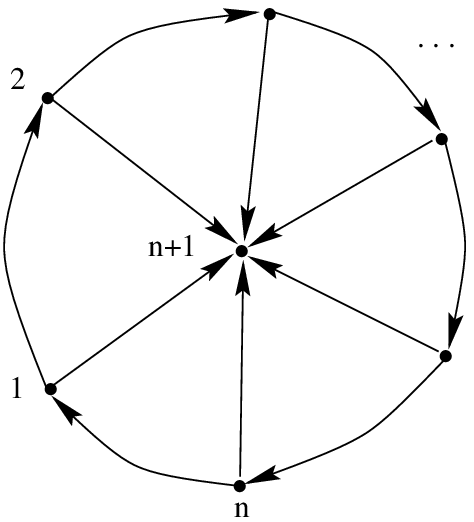}%
\caption{\textsc{The wheel} $\operatorname{Wh}_{n}$}%
\label{-fig-wheel}%
\end{center}
\end{figure}
Arguing as in the proof of Lemma \ref{-estimate-weight}, we conclude that
$w_{\operatorname{Wh}_{n}}\leq4^{n}$, hence $c_{n}^{(1)}\leq2^{n}$.
\end{proof}

By \cite[Th. 8.3.3]{kontsevich}, if $r\in\mathcal{S}(\frak{g})$ is regarded as
a function on $\frak{g}^{\ast}$, then $I_{T}(r)$ is given by the application
of the infinite differential operator associated to the series $S_{1}(x)$ from
Proposition \ref{-tau}. Regarding $r$ as a distribution on $\frak{g}$, we
conclude that
\[
I_{T}(r)=r\tau\text{.}%
\]

We can now proceed with the proof of Theorem \ref{-main-r}.

Let $\Gamma$ be a graph in $A_{n}$ ($V_{\Gamma}=\left\{  1,2,\ldots
,n,\overline{1},\overline{2}\right\}  $). For each configuration in
$C_{n+2,0}$\ (as opposed to $C_{n,2}$ earlier), we can draw $\Gamma$ in
$\mathcal{H}$ and calculate the edge angles. This gives a map
\[
\Phi_{\Gamma}^{\prime}:\overline{C}_{n+2,0}\rightarrow\mathbb{T}^{2n}\text{,}%
\]
and we denote by $\omega_{\Gamma}^{\prime}$ the pullback of the volume form on
the torus.

Next recall that $\overline{C}_{2,0}$ (the \textsc{Eye}) is the compactified
configuration space of two points in $\mathcal{H}$. We fix a path
$\xi:[0,1]\rightarrow\overline{C}_{2,0}$, such that
\begin{align*}
\xi(0)  &  \in C_{2}\text{ (two points coinciding in }\mathcal{H}\text{)},\\
\xi(1)  &  \in C_{0,2}^{+}\text{ (both points lying in }\mathbb{R\subset
}\overline{\mathcal{H}}\text{).}%
\end{align*}%

\begin{figure}
[h]
\begin{center}
\includegraphics[
trim=0.000000in 0.149533in 0.000000in 0.131488in,
height=1.0144in,
width=1.8144in
]%
{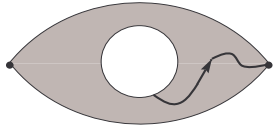}%
\caption{\textsc{ Path in the Eye}}%
\end{center}
\end{figure}

Now for each $n\geq0$ we have a ``forgetting'' map $\overline{C}%
_{n+2,0}\rightarrow\overline{C}_{2,0},$ where we forget the location of first
$n$ points in the configuration. We denote by $Y_{n}\subset$ $\overline
{C}_{n+2,0}$ the preimage of the path $\xi\left(  \lbrack0,1]\right)  $ under
the forgetting map and by $\partial Y_{n}$ its boundary.

The boundary $\partial Y_{n}$ can be written as the following cycle:
\[
\partial Y_{n}=Z_{n}^{0}-Z_{n}^{1}-Z_{n},
\]
where $Z_{n}^{0},Z_{n}^{1}$ are the preimages of the points $\xi(0)$ and
$\xi(1)$ respectively and $Z_{n}$ is the sum (with appropriate signs) of the
components of
\[
\partial\overline{C}_{n+2,0}\cap\text{preimage of }\left\{  \xi
(t):\,0<t<1\right\}  \text{.}%
\]

Since $\omega_{\Gamma}^{\prime}$ is a closed $2n$-form, by Stokes' theorem we
get
\[
\int_{Z_{n}^{0}}\omega_{\Gamma}^{\prime}-\int_{Z_{n}^{1}}\omega_{\Gamma
}^{\prime}=\int_{Z_{n}}\omega_{\Gamma}^{\prime}.
\]

\begin{proof}
[Proof of Theorem \ref{-main-r}]By \cite[8.2.1, 8.2.2]{kontsevich}, for
$r,p\in\mathcal{S}$ we have
\begin{align*}
r\cdot\left[  \sum_{n=0}^{\infty}\frac{1}{n!}\sum_{\Gamma\in A_{n}}%
\partial_{\Gamma}^{p}\left(  \int\limits_{Z_{n}^{1}}\omega_{\Gamma}^{\prime
}\right)  \right]   &  =(r\tau)\star(p\tau)\text{,}\\
r\cdot\left[  \sum_{n=0}^{\infty}\frac{1}{n!}\sum_{\Gamma\in A_{n}}%
\partial_{\Gamma}^{p}\left(  \int\limits_{Z_{n}^{0}}\omega_{\Gamma}^{\prime
}\right)  \right]   &  =(r\ast_{\frak{g}}p)\tau\text{.}%
\end{align*}
Hence for $T_{p}$ as in the statement of Theorem \ref{-main-r}, we get
\[
r\cdot T_{p}=(r\ast_{\frak{g}}p)\tau-(r\tau)\star(p\tau)=r\cdot\left[
\sum_{n=0}^{\infty}\frac{1}{n!}\sum_{\Gamma\in A_{n}}\partial_{\Gamma}%
^{p}\left(  \int\limits_{Z_{n}}\omega_{\Gamma}^{\prime}\right)  \right]
\text{.}%
\]

We now discuss the structure of $Z_{n}$. The relative positions of the
vertices $\overline{1}$ and $\overline{2}$ are determined by the point
$\xi(t)\in C_{2,0}$. The integral $\int\limits_{Z_{n}}\omega_{\Gamma}^{\prime
}$ is the sum of integrals over the top (i.e., $2n$-dimensional) components of
$Z_{n}$, which correspond to the following (degenerate) configurations:

\begin{enumerate}
\item  Two or more points cluster at $\overline{1}$.

\item  Two or more points cluster at $\overline{2}$.

\item  Two or more points cluster somewhere else in $\mathcal{H}$.

\item  One or more points cluster on $\mathbb{R}$.
\end{enumerate}

For components of type 4 $\omega_{\Gamma}^{\prime}=0$ (since the angle for any
edge originating from $\mathbb{R}$ is identically 0). For components of types
1--3 which involve clusters of \emph{three or more} points, $\int
\omega_{\Gamma}^{\prime}=0$ by \cite[Lemma 6.6]{kontsevich}.

Now let $Z$ be a component corresponding to a two-point cluster of type 3. The
corresponding boundary component of $\overline{C}_{n+2,0}$ is $C_{2}\times
C_{n+1,0}$ and we have $Z\simeq C_{2}\times\widetilde{Z}$, where
$\widetilde{Z}$ is the preimage of the path $\xi$ in $C_{n+1,0}$ under the
forgetting map. As in \cite[6.4.1]{kontsevich}, $\int_{Z}$ $\omega_{\Gamma
}^{\prime}$ decomposes as a product of integrals over $C_{2}$ and
$\widetilde{Z}$. The first integral vanishes unless two points in the cluster
are connected by an edge.
\[
\cdots\overset{k}{\leftarrow}\frame{$\overset{\gamma^{kl}}{\bullet}\overset
{l}{\longrightarrow}\overset{\gamma^{ij}}{\bullet}$}\underset{j}{\overset
{i}{\rightrightarrows}}\cdots\text{ .}%
\]

Relabeling the vertices and edges as above, we see that the operator
$\partial_{\Gamma}^{p}$ is a sum of terms with factors of the form
\[
\sum_{l=1}^{d}\gamma^{kl}\frac{\partial\gamma^{ij}}{\partial x_{l}}=\frac
{1}{4}\sum_{l=1}^{d}[x_{k},x_{l}]c_{ij}^{l}=\frac{1}{4}\left[  x_{k,}\left[
x_{i},x_{j}\right]  \right]  \text{.}%
\]

The set of graphs with the above subgraph can be grouped into sets of three
corresponding to a cyclic permutation of the three outgoing edges from the
cluster. Each of these graphs has the same weight $\int_{Z}$ $\omega_{\Gamma
}^{\prime}$. By the Jacobi identity, the sum of the three corresponding
operators $\partial_{\Gamma}^{p}$ vanishes. Hence $\sum_{\Gamma\in A_{n}%
}\partial_{\Gamma}^{p}\left(  \int_{Z}\omega_{\Gamma}^{\prime}\right)  =0$.

Now we consider two-point clusters of type 2.
\[
\ldots\overset{i}{\leftarrow}\frame{$\overset{\gamma^{ij}}{\bullet}\overset
{j}{\longrightarrow}\overset{p}{\bullet}$}\,\cdots
\]
Relabeling the edges and vertices as above, we see that the corresponding
operator is a sum of terms with factors of the form
\[
\sum_{j=1}^{d}(p\cdot y_{j})\left[  x_{i},x_{j}\right]  \text{,}%
\]
where $\left[  x_{i},x_{j}\right]  \in\frak{g}$ acts as a differential
operator. But this expression is equal to $p\cdot\operatorname*{adj}%
\nolimits_{x_{i}}$ (where $\operatorname*{adj}\nolimits_{x_{i}}$ is the
adjoint vector field corresponding to $x_{i}\in\frak{g}$), which vanishes
since $p$ is invariant ($p\in\mathcal{I}$).

This leaves only two-point clusters of type 1.
\[
\ldots\overset{i}{\leftarrow}\frame{$\overset{\gamma^{ij}}{\bullet}\overset
{j}{\longrightarrow}\overset{r}{\bullet}$}\,\cdots
\]
Arguing as above, we conclude that $r\cdot T_{p}$ is a sum of terms containing
factors of the type $r\cdot\operatorname*{adj}\nolimits_{x_{i}}$. Hence
$T_{p}\ $is an infinite sum of operators from $\frak{R}$. By Theorem
\ref{-main-d} and Proposition \ref{-tau}
\[
T_{p}=\partial_{p}\tau-\tau\partial_{p\tau}^{\star}%
\]
belongs to $\frak{D}$. Therefore, $T_{p}\in\frak{R}$, as was to be shown.
\end{proof}

\textbf{Remark.} In the above discussion we did not explicitly take into
consideration the possibility of edges originating outside the cluster and
terminating at the vertices of the cluster. However, we can group together
graphs which differ only by a permutation of the incident vertices in the
cluster, and observe that these graphs have the same weight. The result now
follows from the Leibniz formula. \label{-tang-for-lie}

\section{Kashiwara-Vergne conjecture and local solvability}

In this section we prove the Kashiwara-Vergne conjecture (Theorem
\ref{-kv-conjecture}) and discuss its consequences.

We start by observing, that for $x^{\prime}$ and $x^{\prime\prime}$ in
$\frak{g}$, we have
\[
\left[  x^{\prime},x^{\prime\prime}\right]  =x^{\prime}\star x^{\prime\prime
}-x^{\prime\prime}\star x^{\prime}.
\]
Then, by the universal property of $\mathcal{U}=$ $\mathcal{U}(\frak{g})$ it
follows that there is a unique algebra homomorphism between $\mathcal{U}$ and
$(\mathcal{S},\star)$. This morphism (denoted $I_{alg}$ in \cite[8.3.1]%
{kontsevich}) is invertible, and we shall denote its inverse by $\kappa$.

Then, as established in \cite[8.3.4]{kontsevich}, we have
\[
\kappa(p\tau)=\exp_{\ast}(pq)=\eta(p)\ \text{for }p\in\mathcal{S}\text{.}%
\]

\begin{proof}
[Proof of Theorem \ref{-kv-conjecture}]In view of Theorem \ref{-main-r}, it
suffices to establish the following identity (in $\frak{D}$)
\[
D_{p}=T_{p}\tau^{-1}q\text{.}%
\]
Since both sides are germs of analytic differential operators, it is enough to
verify that
\[
r\cdot D_{p}=\left(  r\cdot T_{p}\right)  \tau^{-1}q
\]
for all $r\in\mathcal{S}$.

Then
\begin{align*}
r\cdot D_{p}  &  =\exp^{\ast}\left[  \eta(r\ast_{\frak{g}}p)-\eta(r)\ast
_{G}\eta(p)\right] \\
&  =(r\ast_{\frak{g}}p)q-\exp^{\ast}(\kappa(r\tau)\ast_{G}\kappa(p\tau))\\
&  =(r\ast_{\frak{g}}p)q-\exp^{\ast}(\kappa(r\tau\star p\tau))\\
&  =(r\ast_{\frak{g}}p)q-\exp^{\ast}\left(  \eta\left(  (r\tau\star p\tau
)\tau^{-1}\right)  \right) \\
&  =(r\ast_{\frak{g}}p)q-(r\tau\star p\tau)\tau^{-1}q\\
&  =\left[  (r\ast_{\frak{g}}p)\tau-(r\tau\star p\tau)\right]  \tau^{-1}q\\
&  =\left(  r\cdot T_{p}\right)  \tau^{-1}q\text{.}%
\end{align*}
\end{proof}

Theorems \ref{-module-iso} and \ref{-eigendist} follow immediately from the
Theorem above, and Theorem \ref{-fund-sol} requires a short argument, which we
give below.

\begin{proof}
[Proof of Theorem \ref{-fund-sol}]A bi-invariant differential operator on $G$
is precisely an element of the center $\mathcal{Z}$ of $\mathcal{U}$, hence of
the form $\eta(p)$ for some $p\in\mathcal{I}$. Let $\delta_{0}$ and
$\delta_{1}$ be the delta-distributions supported at $0\in\frak{g}$ and $1\in
G$ respectively. Then if $p\in\mathcal{I}\mathbf{\setminus}\{0\}$, by
\cite{rais-nilpotent} (see also \cite{rais-review}), there exists a germ
$P\in\mathbf{I}$, such that
\[
P\ast_{\frak{g}}p=\delta_{0}\text{.}%
\]
Hence by Theorem \ref{-module-iso}, we get
\[
\eta(P)\ast_{G}\eta(p)=\eta(\delta_{0})=\delta_{1}\text{,}%
\]
and we see that $\eta(P)$ is a local fundamental solution for $\eta(p)$.\label{-sec-kv}
\end{proof}

\section{Results for Lie supergroups}

\subsection{Preliminaries on supermanifolds}

We briefly discuss basic facts about supermanifolds. The study of
supermanifolds etc. was initiated by Berezin (see \cite{berezin},
\cite{manin}) and independently by Kostant (\cite{kostant}). Subsequently a
slightly different (though essentially equivalent) treatment appeared in the
physics literature \cite{super}. We will follow the approach of \cite{kostant}%
. The reader familiar with the physics point of view should have no trouble
translating everything into that language.

A supermanifold ($\mathbb{Z}_{2}$-graded manifold) of dimension $(d_{0}%
,d_{1})$ is a triple $(X,\mathcal{A},\pi),$ where $X$ is a $d_{0}$-dimensional
smooth manifold, $\mathcal{A}$ is a sheaf of $\mathbb{Z}_{2}$-graded
commutative algebras and $\pi$ is a sheaf map from $\mathcal{A}$ to the sheaf
$\mathcal{C}^{\infty}(X)$ of smooth functions on $X$. Moreover, we require
that every open subset of $X$ can be covered by ``$\mathcal{A}$-splitting open
sets $U$ of odd dimension $d_{1}$''. This means that we can choose
(non-canonical) subalgebras $C(U)$ and $D(U)$ of $\mathcal{A}(U)$ such that
\begin{gather*}
\pi|_{C(U)}:C(U)\longrightarrow C^{\infty}(U)\text{ is an algebra
isomorphism,}\\
D(U)\text{ is isomorphic to the exterior algebra in }d_{1}\text{ variables,}\\
C(U)\otimes D(U)\simeq\mathcal{A}(U)\text{ (algebra isomorphism).}%
\end{gather*}

$\mathcal{A}(X)$ and its dual $\mathcal{A}(X)^{\prime}=\operatorname{Hom}%
_{\mathbb{R}}(\mathcal{A}(X),\mathbb{R})$ are super (i.e., $\mathbb{Z}_{2}%
$-graded) vector spaces. A ``pure'' linear functional $v\in\mathcal{A}%
(X)^{\prime}$ is called a differentiation at $x\in X$ if%
\[
\left\langle v,fg\right\rangle =\left\langle v,f\right\rangle \widetilde
{g}(x)+(-1)^{p(v)p(f)}\widetilde{f}(x)\left\langle v,g\right\rangle ,
\]
for any pure $f$ and $g$ in $\mathcal{A}(X)$. Here we write $\widetilde{f}$
for $\pi(f)$ and $p$ for the parity ($0$ or $1$). An arbitrary $v\in
\mathcal{A}(X)^{\prime}$ is called a differentiation if its homogeneous
components are differentiations, and we write $T_{x}(X,\mathcal{A})$ for the
super vector space of all differentiations at $x$. It is easy to see that
$\dim\left[  T_{x}(X,\mathcal{A})\right]  _{i}=d_{i}$, $i=0,1$. Moreover, we
have a natural surjection $\pi_{x}:$ $T_{x}(X,\mathcal{A})\longrightarrow
T_{x}(X)$ (ordinary tangent space of $X$ at $x$) with $\ker\pi_{x}=\left[
T_{x}(X,\mathcal{A})\right]  _{1}$.

For each $k=0,1,\ldots$ the sheaf $\operatorname*{Diff}_{k}\mathcal{A}$ of
differential operators of order $\leq k$ on the supermanifold $(X,\mathcal{A)}%
$ is defined in \cite[2.10]{kostant} as follows:%
\begin{align*}
\left(  \operatorname*{Diff}\nolimits_{0}\mathcal{A}\right)  (U)  &
=\mathcal{A}(U)\text{, regarded as left multiplication operators;}\\
\left(  \operatorname*{Diff}\nolimits_{k}\mathcal{A}\right)  (U)  &  =\left\{
\partial\in\operatorname*{End}\mathcal{A}(U):\left[  \partial,f\right]
\in\left(  \operatorname*{Diff}\nolimits_{k-1}\mathcal{A}\right)  (U)\text{
for all }f\in\mathcal{A}(U)\right\}  \text{.}%
\end{align*}
Here the bracket $\left[  \ ,\ \right]  $ denotes the graded commutator in
$\operatorname*{End}\mathcal{A}(U)$.

For $\partial\in\left(  \operatorname*{Diff}\nolimits_{k}\mathcal{A}\right)
(X)$ we can define its \emph{body} as the unique differential operator
$\widetilde{\partial}$ on $X$ satisfying the conditions%
\begin{align*}
\widetilde{\partial}\cdot1_{C^{\infty}(X)}  &  =\widetilde{\partial
\cdot1_{\mathcal{A}}}\text{ and}\\
\left[  \widetilde{\partial},\widetilde{f}\right]   &  =\widetilde{\left[
\partial,f\right]  }\text{ \ for all }f\in\mathcal{A}(X)\text{.}%
\end{align*}

This allows us to introduce a locally convex topology on $\mathcal{A}(X)$, as
in \cite[2.9]{kostant}. The space of distributions $\operatorname*{Dist}%
(X,\mathcal{A})$ can then be defined as the topological dual of $\mathcal{A}%
(X)$. The support of a function $f\in\mathcal{A}(X)$ is simply the support of
$\pi(f)\in C^{\infty}(X)$, and the support of a distribution $T\in$
$\operatorname*{Dist}(X,\mathcal{A})$ is the smallest closed set
$\operatorname*{supp}T$ such that $T$ vanishes on all $f\in\mathcal{A}(X)$
with $\operatorname*{supp}f\cap\operatorname*{supp}T=\emptyset$.

\subsection{Preliminaries on Lie supergroups}

Following \cite[2.11]{kostant}, we define $\mathcal{A}(X)^{\ast}%
\subset\mathcal{A}(X)^{\prime}\ $as the subset of those linear functionals on
$\mathcal{A}(X)$ whose kernel contains an ideal of finite codimension. It is
easy to check that $\mathcal{A}(X)^{\ast}$ is a graded co-commutative
coalgebra. We say that the supermanifold $(G,\mathcal{A},\pi)$ is a Lie
supergroup (graded Lie group in \cite{kostant}), if $\mathcal{A}(G)^{\ast}$
has the structure of a Hopf algebra extending the natural coalgebra structure.

This implies that $G$ is a Lie group and that the super tangent space
$\frak{g}=T_{e}(G,\mathcal{A})$ at the identity $e\in G$ has the structure of
a Lie superalgebra, i.e.,%
\[
\frak{g}=\frak{g}_{0}+\frak{g}_{1}\text{,}%
\]
where $\frak{g}_{0}$ is the Lie algebra of $G$. Moreover, we have a
representation $\pi$ of $G$ on $\frak{g}_{1}$, whose differential is the
adjoint action of $\frak{g}_{0}$ on $\frak{g}_{1}$. We have used the same
letter $\pi$ to emphasize the fact that the supergroup can be reconstructed
from the representation. The Hopf algebra $\mathcal{A}(G)^{\ast}$ can be
obtained as the smash product of the group algebra $\mathbb{R}G$ and the
universal enveloping algebra $\mathcal{U}(\frak{g})$, and then $\mathcal{A}%
(G)$ can be recovered as in \cite[3.7]{kostant}.

The action of a Lie supergroup $(G,\mathcal{A})$ on a supermanifold is defined
in \cite[3.9]{kostant}. Each Lie supergroup $(G,\mathcal{A})$ has an ``adjoint
action'' on itself and on its Lie superalgebra $\frak{g}$ (regarded as a
supermanifold). Thus for each $a\in\frak{g}\frak{\ }$we obtain derivations of
the algebras $\mathcal{A}(G)$ and $\mathcal{A}(\frak{g})$, which will both be
denoted by $\operatorname*{adj}\nolimits_{a}$.

As before, we regard the (graded) symmetric algebra $\mathcal{S}%
=\mathcal{S}(\frak{g})$ as the convolution algebra of distributions on the
super vector space $\frak{g}$ with support at $0\in\frak{g}$ and denote by
$\mathcal{I}$ the subalgebra of $\operatorname*{adj}$--invariant distributions.

We regard the universal enveloping algebra $\mathcal{U}=\mathcal{U}(\frak{g})$
as the convolution algebra of distributions on the Lie supergroup
$(G,\mathcal{A},\pi)$ with support at $e\in G$. The center $\mathcal{Z}$ of
$\mathcal{U}$ can then be regarded as the subalgebra of $\operatorname*{adj}%
$--invariant distributions.

The spaces of germs $\mathbf{S}$, $\mathbf{I}$, $\mathbf{U}$, $\mathbf{Z}%
\ $(see Introduction) can be defined analogously in the super-setting.

The exponential map from $\frak{g}$ to $(G,\mathcal{A})$ is best described in
terms of the corresponding ``pushforward'' map
\[
\exp_{\ast}:\mathcal{A}(\frak{g})^{\ast}\longrightarrow\mathcal{A}(G)^{\ast}%
\]
which is the usual exponential from $\mathbb{R}\frak{g}$ to the group algebra
$\mathbb{R}G$, and the symmetrization map from $\mathcal{S}$ to $\mathcal{U}$.

At the germ level $\exp_{\ast}$ is an isomorphism from $\mathbf{S}$ to
$\mathbf{U}$ (and $\mathbf{I}$ to $\mathbf{Z}$). Hence it has a well-defined
inverse $\exp^{\ast}$.

As in the Introduction, we denote by $\frak{D}$ the algebra of germs at $0$ of
differential operators on $\frak{g}$ with analytic coefficients, and by
$\frak{R}$ the right ideal of $\frak{D}$ generated by the germs of the
differential operators $\operatorname*{adj}{}_{a}$, $a\in\frak{g}$.

\subsection{Statements and proofs for Lie supergroups}

All results stated in Sections \ref{-sec-intro} to \ref{-sec-kv} hold with
minor changes in the statements and proofs. We discuss the less obvious
modifications below.

Let $\frak{g}$ be a finite-dimensional real Lie superalgebra. Then its dual
\[
\frak{g}^{\ast}=\operatorname{Hom}_{\mathbb{R}}(\frak{g},\mathbb{R}%
)=\frak{g}_{0}^{\ast}+\frak{g}_{1}^{\ast}%
\]
is a Poisson supermanifold (cf. \cite[5.2]{kostant}). The Poisson bracket
$\gamma$ defined as the unique biderivation of $\mathcal{A}(\frak{g}^{\ast})$
satisfying%
\[
\gamma(f_{1},f_{2})=\left\{
\begin{array}
[c]{ll}%
0 & \text{if }f_{1}\text{ or }f_{2}\text{ is constant}\\
\frac{1}{2}\left[  f_{1},f_{2}\right]  & \text{if }f_{1},f_{2}\in\frak{g}%
\end{array}
\right.  \text{.}%
\]
For convenience we fix bases $x_{1},x_{2},\ldots,x_{d_{0}}$ of $\frak{g}_{0}$
and $x_{d_{0}+1},\ldots,x_{d_{0}+d_{1}}$ of $\frak{g}_{1}.$ Then in
coordinates we have
\[
\gamma(f_{1},f_{2})=\frac{1}{2}\sum_{i,j}[x_{i},x_{j}]\frac{\partial f_{1}%
}{\partial x_{i}}\frac{\partial f_{2}}{\partial x_{j}}\text{.}%
\]
For an odd variable $x_{i}$ the vector field $\frac{\partial}{\partial x_{i}}$
is odd ; nevertheless, $\gamma$ is an \emph{even} bidifferential operator and
preserves the grading on $\mathcal{A}(\frak{g}^{\ast})$.

Following \cite{kontsevich}, we can define the $\star$-product of $f_{1}%
,f_{2}\in\mathcal{A}(\frak{g}^{\ast})$ by the formula%
\[
f_{1}\star f_{2}=\sum_{n=0}^{\infty}\frac{\hbar^{n}}{n!}\sum_{\Gamma\in G_{n}%
}w_{\Gamma}B_{\Gamma}(f_{1},f_{2})\text{.}%
\]
The bidifferential operators $B_{\Gamma}$ for $\Gamma\in G_{n}$ are given by
the formula (\ref{=ouch}) up to a sign which can be calculated in the
following way (illustrated by the example below).

We consider the same graph $\Gamma_{0}\in G_{2}$ as in Section
\ref{ss-weights}, i.e.%

\[
\text{\ }%
\begin{array}
[c]{ccc}%
\overset{1}{\bullet} & \rightarrow & \overset{2}{\bullet}\\
\downarrow & \swarrow & \downarrow\\
\underset{\overline{1}}{\bullet} &  & \underset{\overline{2}}{\bullet}%
\end{array}
\text{\ or (relabeling)\ \ \ }%
\begin{array}
[c]{ccc}%
\overset{\gamma^{i_{1}j_{1}}}{\bullet} & \overset{_{j_{1}}}{\rightarrow} &
\overset{\gamma^{i_{2}j_{2}}}{\bullet}\\
\,\downarrow^{_{i_{1}}} & \overset{i_{2}}{\swarrow} & \downarrow^{_{j_{2}}}\\
\underset{f_{1}}{\bullet} &  & \underset{f_{2}}{\bullet}%
\end{array}
.
\]
We first form the expression
\begin{equation}
\sum_{1\leq i_{1},j_{1},i_{2},j_{2}\leq d}(\partial_{i_{1}}\gamma^{i_{1}j_{1}%
}\partial_{j_{1}})(\partial_{i_{2}}\gamma^{i_{2}j_{2}}\partial_{j_{2}}%
)f_{1}\,f_{2}\text{,} \label{=super-mnemonic}%
\end{equation}
where the bivector fields $(\gamma)$ corresponding to the vertices $1$ and $2$
are all written on the left (since $\gamma$ is even, the order is irrelevant).
Then we use the Koszul rule of signs to move the derivatives to the
appropriate positions (as given by the graph), i.e. $\partial_{i_{1}}$ and
$\partial_{i_{2}}$ should be moved to the left of $f_{1}$, $\partial_{j_{2}}$
to $f_{2}$ and $\partial_{j_{1}}$ to $\gamma^{i_{2}j_{2}}$. We obtain the
bidifferential operator, as in the formula (\ref{=ouch-example}) before
\begin{equation}
\sum_{1\leq i_{1},j_{1},i_{2},j_{2}\leq d}\pm\gamma^{i_{1}j_{1}}%
\partial_{j_{1}}(\gamma^{i_{2}j_{2}})\,\,\partial_{i_{1}}\partial_{i_{2}%
}(f_{1})\,\,\partial_{j_{2}}(f_{2})\text{.} \label{=super-term}%
\end{equation}
Here the signs of the summands are given by $-1$ to the power%
\[
p(j_{2})p(f_{1})+p(i_{2})\left(  p(i_{2})+p(j_{2})\right)  +p(i_{1})\left(
p(i_{1})+p(i_{2})+p(j_{2})\right)  \text{, }%
\]
where $p(i)$ is the parity of $x_{i}\in\frak{g}$ and $p(f)$ is the parity of
the $f\in\mathcal{A}(\frak{g}^{\ast})$. Note that in the above formula
(\ref{=super-term}) we did not move $\partial_{i_{1}}$ past $\partial_{i_{2}}%
$, thus preserving the same order as in (\ref{=super-mnemonic}).

It is easy to check that Kontsevich's argument for the associativity of the
$\star$-product remains valid for $\mathcal{A}(\frak{g}^{\ast})$.

Proceeding as in Section \ref{-star-for-lie}, we obtain

\begin{theorem}
Given $p\in\mathcal{S}$ of order $l$, there is a unique element $\partial
_{p}^{\star}$ of order $l$ in $\frak{D}$ such that for all $r\in\mathcal{S}$%
\[
r\star p=r\cdot\partial_{p}^{\star}\text{ .}%
\]
\end{theorem}

The analyticity of the coefficients of $\partial_{p}^{\star}$ is verified
exactly as before.

Let us write $\operatorname{str}$ for the supertrace of an endomorphism of a
graded vector space \cite{manin}. The expression $\operatorname{str}\left[
\left(  \operatorname*{ad}x\right)  ^{2k}\right]  $ defines an even element in
$\mathcal{A}(\frak{g})$. With $c_{2k}^{(1)}$ as in (\ref{=s1}), we can verify
(as in Section \ref{-tang-for-lie}) that
\[
S_{1}(x)=\exp\left(  \sum_{k=1}^{\infty}c_{2k}^{(1)}\operatorname{str}\left[
\left(  \operatorname*{ad}x\right)  ^{2k}\right]  \right)
\]
converges to an even superfunction $\tau(x)$ in $\mathcal{A}(\frak{g})$ which
is analytic at $0\in\frak{g}$.

Proceeding as in Section \ref{-tang-for-lie}, we obtain

\begin{theorem}
For $p\in I$, the operator
\[
T_{p}=\partial_{p}\tau-\tau\partial_{p\tau}^{\star}%
\]
belongs to $\frak{R}$.\label{-supermain-r}
\end{theorem}

Just as for $\tau(x)$, the expression%
\[
q(x)=\exp\left(  \sum_{k=1}^{\infty}\frac{B_{2k}}{4k\,(2k)!}\operatorname{str}%
\left[  \left(  \operatorname*{ad}x\right)  ^{2k}\right]  \right)  \text{.}%
\]
defines an element in $\mathcal{A}(\frak{g}),$ which is analytic on $\frak{g}$.

The map $\eta:\mathbf{S\longrightarrow U}$ is defined (as before) by
\[
\eta(p)=\exp_{\ast}(pq)\text{,}%
\]
and the restriction of $\eta$ to $\mathcal{I}$ is an \emph{algebra}
isomorphism from $\mathcal{I}$ to $\mathcal{Z}$ (this extension of Duflo's
isomorphism to the case of Lie superalgebras is due to Kontsevich).

Theorems \ref{-kv-conjecture}, \ref{-module-iso} and \ref{-eigendist} for any
real Lie superalgebra $\frak{g}$ now follow from Theorem \ref{-supermain-r},
as in Section \ref{-sec-kv}.

The statement of Theorem \ref{-fund-sol} should be modified in the following manner:

\begin{theorem}
Any bi-invariant differential operator on a real Lie supergroup $G$, with a
nonzero body, admits a (local) fundamental solution and hence is locally solvable.
\end{theorem}

\begin{proof}
If a bi-invariant differential operator on a real Lie supergroup $G$ has a
nonzero body, it can be represented as $\eta(p)$, where $p\in\mathcal{I}$ and
$p$ has a nonzero body. Then we can find a germ $P\in\mathbf{I}$, such that
\[
P\ast_{\frak{g}}p=\delta_{0}\text{.}%
\]
Local solvability of $\eta(p)$ now follows as in the proof of Theorem
\ref{-fund-sol}.
\end{proof}

\end{document}